# M-ESTIMATION OF LINEAR MODELS WITH DEPENDENT ERRORS[1]

By Wei Biao Wu

*University of Chicago*

We study asymptotic properties of $M$-estimates of regression parameters in linear models in which errors are dependent. Weak and strong Bahadur representations of the $M$-estimates are derived and a central limit theorem is established. The results are applied to linear models with errors being short-range dependent linear processes, heavy-tailed linear processes and some widely used nonlinear time series.

**1. Introduction.** Consider the linear model

(1) $$y_i = \mathbf{x}_i'\beta + e_i, \qquad 1 \le i \le n,$$

where $\beta$ is a $p \times 1$ unknown regression coefficient vector, $\mathbf{x}_i = (x_{i1}, \ldots, x_{ip})'$ are $p \times 1$ known (nonstochastic) design vectors and $e_i$ are errors. We estimate the unknown parameter vector $\beta$ by minimizing

(2) $$\sum_{i=1}^{n} \rho(y_i - \mathbf{x}_i'\beta),$$

where $\rho$ is a convex function. Important examples include Huber's estimate with $\rho(x) = (x^2 \mathbf{1}_{|x| \le c})/2 + (c|x| - c^2/2)\mathbf{1}_{|x| > c}$, $c > 0$, the $\mathcal{L}^q$ regression estimate with $\rho(x) = |x|^q$, $1 \le q \le 2$, and regression quantiles with $\rho(x) = \rho_\alpha(x) = \alpha x^+ + (1-\alpha)(-x)^+$, $0 < \alpha < 1$, where $x^+ = \max(x, 0)$. In particular, if $q = 1$ or $\alpha = 1/2$, then the minimizer of (2) is called the least absolute deviation (LAD) estimate. See [2] and [68] for $\mathcal{L}^q$ regression estimates and [35] for regression quantiles. See also [34] for an excellent account of quantile regression.

Received November 2005; revised May 2006.
[1]Supported by NSF Grant DMS-04-78704.
*AMS 2000 subject classifications.* Primary 62J05; secondary 60F05.
*Key words and phrases.* Robust estimation, linear model, dependence, nonlinear time series.

Let $\hat{\beta}_n$ be the minimizer of (2) and let $\beta_0$ be the true parameter. There is a substantial amount of work concerning asymptotic properties of $\hat{\beta}_n - \beta_0$ for various forms of $\rho$ (not necessarily convex); see, for example, [2, 3, 4, 7, 8, 10, 12, 26, 30, 32, 49, 57, 66, 67] and [69] among others. Deep results such as Bahadur representations have also been obtained. However, in the majority of the previous work it is assumed that the errors $e_i$ are independent. The asymptotic problem of $M$-estimation of linear models with dependent errors is practically important, however theoretically challenging. Huber [29, 30] commented that the assumption of independence is a serious restriction. See also [24].

In this paper we shall relax the independence assumption in the classical $M$-estimation theory so that a very general class of dependent errors is allowed. Specifically, we shall establish a Bahadur representation and a central limit theorem for $\hat{\beta}_n - \beta_0$ for the linear model (1) with the errors $(e_i)$ being stationary causal processes. In the early literature very restrictive assumptions were imposed on the error process $(e_i)$. Typical examples are strongly mixing processes of various types. See [13, 36] and [41] among others for strong mixing processes and [45] for $\varphi$-mixing processes. Berlinet, Liese and Vajda [9] obtained consistency of $M$-estimators for regression models with strong mixing errors. Gastwirth and Rubin [20] considered the behavior of $L$-estimators of strong mixing Gaussian processes and first-order autoregressive processes with double exponential marginals. It is generally not easy to verify strong mixing conditions. For example, for linear processes to be strong mixing, very restrictive conditions are needed on the decay rate of the coefficients [17, 22, 40, 60]. Portnoy [43, 44] and Lee and Martin [38] investigated the effect of dependence on robust location estimators by assuming that the errors are autoregressive moving average processes with finite orders.

To the best of our knowledge, it seems that the problem of Bahadur representations has been rarely studied for $M$-estimates of linear models with nonstrong mixing errors. The Bahadur-type representations provide significant insight into the asymptotic behavior of an estimator by approximating it by a linear form. For sample quantiles it has been investigated by Hesse [27], Babu and Singh [5] and Wu [62], among others. Babu [4] considered LAD estimators for linear models with strong mixing errors.

For the errors $(e_i)$ we confine ourselves to stationary causal processes. Namely, let

$$(3) \qquad e_i = G(\ldots, \varepsilon_{i-1}, \varepsilon_i),$$

where $\varepsilon_k$, $k \in \mathbb{Z}$, are independent and identically distributed (i.i.d.) random variables and $G$ is a measurable function such that $e_i$ is a proper random variable. Here $\mathbb{Z}$ is the set of integers. The framework (3) is a natural paradigm for nonlinear time series models and it represents a huge class of stationary processes which appear frequently in practice. As in [46, 53, 59] and [63], (3) can be interpreted as a physical system with the innovations $\varepsilon_i$ being the inputs that drive the system, $G$ being a filter and $e_i$ being the output. This interpretation leads to our dependence measures. The Wiener conjecture states that every stationary and ergodic process $(e_i)$ can be expressed in the form of (3); see [33, 50, 51] and [55], page 204.

Let the shift process $\mathcal{F}_k = (\ldots, \varepsilon_{k-1}, \varepsilon_k)$. For $i \in \mathbb{N}$ let $F_i(u|\mathcal{F}_0) = \mathbb{P}(e_i \leq u|\mathcal{F}_0)$ [resp., $f_i(u|\mathcal{F}_0)$] be the conditional distribution (resp., density) function of $e_i$ at $u$ given $\mathcal{F}_0$ and let $f$ be the marginal density of $e_i$. Let $l \geq 0$. For a function $g$, write $g \in \mathcal{C}^l$ if $g$ has derivatives up to $l$th-order and $g^{(l)}$ is continuous. Denote by $f_i^{(l)}(u|\mathcal{F}_0) = \partial^l f_i(u|\mathcal{F}_0)/\partial u^l$ the $l$th-order derivative if it exists. Let $(\varepsilon_i')$ be an i.i.d. copy of $(\varepsilon_i)$, $\mathcal{F}_k^* = (\ldots, \varepsilon_{-1}, \varepsilon_0', \varepsilon_1, \ldots, \varepsilon_k)$ and $e_k^* = G(\mathcal{F}_k^*)$. Then $\mathcal{F}_k^*$ is a coupled version of $\mathcal{F}_k$ with $\varepsilon_0$ replaced by $\varepsilon_0'$, $\mathcal{F}_j^* = \mathcal{F}_j$, $j < 0$, and $e_k$ and $e_k^*$ are identically distributed. Our short-range dependence (SRD) conditions suggest that a certain distance between the two predictive distributions $[e_i|\mathcal{F}_0]$ and $[e_i^*|\mathcal{F}_0^*]$ is summable over $i \geq 1$. Since those conditions are directly related to the data-generating mechanism of $(e_i)$, they are often easily verified; see applications in Section 3.

The paper is structured as follows. Section 2 presents our main results on Bahadur representations and central limit theorems for $\hat{\beta}_n - \beta_0$. Section 3 contains applications to linear models with errors being short-range dependent linear processes, heavy-tailed linear processes where $M$-estimation is particularly relevant, and some widely used nonlinear time series. Proofs are given in Section 4.

**2. Main results.** Without loss of generality, assume throughout the paper that the true parameter $\beta_0 = 0$. We first introduce some notation. Let $\lceil a \rceil = \min\{k \in \mathbb{Z} : k \geq a\}$ and $\lfloor a \rfloor = \max\{k \in \mathbb{Z} : k \leq a\}$, $a \in \mathbb{R}$, be the usual ceiling and floor functions. For a $p$-dimensional vector $\mathbf{v} = (v_1, \ldots, v_p)$ let $|\mathbf{v}| = (\sum_{i=1}^p v_i^2)^{1/2}$. A random vector $V$ is said to be in $\mathcal{L}^q$, $q > 0$, if $\mathbb{E}(|V|^q) < \infty$. In this case write $\|V\|_q = [\mathbb{E}(|V|^q)]^{1/q}$ and $\|V\| = \|V\|_2$. Let the covariance matrix of a $p$-dimensional column random vector $V$ be $\operatorname{var}(V) = \mathbb{E}(VV') - \mathbb{E}(V)\mathbb{E}(V')$. Define projection operators $\mathcal{P}_k$, $k \in \mathbb{Z}$, by $\mathcal{P}_k V = \mathbb{E}(V|\mathcal{F}_k) - \mathbb{E}(V|\mathcal{F}_{k-1})$, $V \in \mathcal{L}^1$. The symbol $C$ denotes a generic constant which may vary from place to place. For a sequence of random variables $(\eta_n)$ and a positive sequence $(d_n)$, write $\eta_n = o_{\text{a.s.}}(d_n)$ if $\eta_n/d_n$ converges to 0 almost surely and $\eta_n = O_{\text{a.s.}}(d_n)$ if $\eta_n/d_n$ is almost surely bounded. We can similarly define the relations $o_\mathbb{P}$ and $O_\mathbb{P}$. Let $N(\mu, \Sigma)$ denote a multivariate normal distribution with mean vector $\mu$ and covariance matrix $\Sigma$.

Let the model matrix $\mathbf{X}_n = (\mathbf{x}_1, \ldots, \mathbf{x}_n)'$ and $\Sigma_n = \mathbf{X}_n'\mathbf{X}_n$. Assume that $\Sigma_n$ is nonsingular for large $n$. It is convenient to consider the rescaled model

$$y_i = \mathbf{z}_{i,n}'\theta + e_i, \tag{4}$$

where $\mathbf{z}_{i,n} = \Sigma_n^{-1/2}\mathbf{x}_i$ and $\theta = \theta_n = \Sigma_n^{1/2}\beta$. Studying the asymptotic behavior of $\hat{\beta}_n$ is equivalent to studying that of $\hat{\theta}_n = \Sigma_n^{1/2}\hat{\beta}_n$, which is a minimizer of $\sum_{i=1}^{n} \rho(e_i - \mathbf{z}'_{i,n}\theta)$. If there are multiple minimizers, we just choose any such minimizer. Observe that $\sum_{i=1}^{n} \mathbf{z}_{i,n}\mathbf{z}'_{i,n} = \mathrm{Id}_p$, the $p \times p$ identity matrix. For $q > 0$ define

$$(5) \qquad \zeta_n(q) = \sum_{i=1}^{n} |\mathbf{z}_{i,n}|^q \quad \text{and} \quad \xi_n(q) = \sum_{i=1}^{n} |\mathbf{x}_i|^q.$$

Assume that $\rho$ has derivative $\psi$. Define the $k$th-step-ahead predicted function

$$(6) \qquad \psi_k(t; \mathcal{F}_0) = \mathbb{E}[\psi(e_k + t)|\mathcal{F}_0], \qquad k \geq 0.$$

The function $\psi_k(\cdot;\cdot)$ plays an important role in the study of the asymptotic behavior of $\hat{\beta}_n$. We now list some regularity conditions on $\rho$, $\mathbf{x}_i$ and the errors $e_i$:

(A1) $\rho$ is a convex function, $\mathbb{E}[\psi(e_1)] = 0$ and $\|\psi(e_1)\|^2 > 0$.
(A2) $\varphi(t) := \mathbb{E}[\psi(e_1 + t)]$ has a strictly positive derivative at $t = 0$.
(A3) $m(t) := \|\psi(e_1 + t) - \psi(e_1)\|$ is continuous at $t = 0$.
(A4) $r_n := \max_{i \leq n} |\mathbf{z}_{i,n}| = \max_{i \leq n}(\mathbf{x}'_i\Sigma_n^{-1}\mathbf{x}_i)^{1/2} = o(1)$.

Conditions (A1)–(A4) are standard and they are often imposed in the $M$-estimation theory of linear models with independent errors; see, for example, [7]. In (A1), the error process $(e_i)$ itself is allowed to have infinite variance, which is actually one of the primary reasons for robust estimation. Section 3 contains an application to linear models with dependent heavy-tailed errors. Under (A2), $\theta$ is estimable or separable. Condition (A3) is very mild. Note that $\psi$ is nondecreasing and it has countably many discontinuity points. If $e_i$ has a continuous distribution function and $\|\psi(e_1 + t_0)\| + \|\psi(e_1 - t_0)\| < \infty$ for some $t_0 > 0$, then $\lim_{t \to 0} \psi(e_1 + t) = \psi(e_1)$ almost surely and (A3) follows from the Lebesgue dominated convergence theorem.

The uniform asymptotic negligibility condition (A4) is basically the Lindeberg–Feller-type condition. With (A4), the diagonal elements of the hat matrix $\mathbf{X}_n\Sigma_n^{-1}\mathbf{X}'_n$ are uniformly negligible. Let $\mathbf{x}_{i_1}, \ldots, \mathbf{x}_{i_p}$ be linearly independent, $1 \leq i_1 < \cdots < i_p$, and $Q = (\mathbf{x}_{i_1}, \ldots, \mathbf{x}_{i_p})$. Then $Q$ is nonsingular, $Q'\Sigma_n^{-1}Q \to 0$ and consequently $\Sigma_n^{-1} \to 0$. The latter implies that the minimum eigenvalue of $\Sigma_n$ diverges to $\infty$ and it is a classical condition for weak consistency of the least squares estimators [18]. For the regression model (1) with i.i.d. errors $e_i$ having mean 0 and finite variance, (A4) is necessary and sufficient for the least squares estimator $\Sigma_n^{-1}\mathbf{X}'_n(y_1, \ldots, y_n)'$ to be asymptotically normal (see [30], Section 7.2, and [21]).

Besides the classical conditions (A1)–(A4), to obtain asymptotic properties of $\hat{\beta}_n$ and $\hat{\theta}_n$ we certainly need appropriate dependence conditions [cf. (7) and (14)]. They are expressed in terms of $\psi_k(\cdot; \mathcal{F}_0)$. Recall $\mathcal{F}_k^* = (\ldots, \varepsilon_{-1}, \varepsilon_0', \varepsilon_1, \ldots, \varepsilon_k)$ and $e_k^* = G(\mathcal{F}_k^*)$.

2.1. *Asymptotic normality.* Theorem 1 asserts that $\hat{\theta}_n$ can be approximated by the linear form $T_n = \sum_{i=1}^n \psi(e_i) \mathbf{z}_{i,n}$ with an $o_\mathbb{P}(1)$ error. Due to the linearity it is easier to deal with $T_n$, which is asymptotically normal under proper conditions (cf. Lemma 2).

THEOREM 1. *Assume* (A1)–(A4) *and, for some* $\epsilon_0 > 0$,

$$(7) \qquad \sum_{i=0}^\infty \sup_{|\epsilon| \leq \epsilon_0} \|\mathbb{E}[\psi(e_i + \epsilon)|\mathcal{F}_0] - \mathbb{E}[\psi(e_i^* + \epsilon)|\mathcal{F}_0^*]\| < \infty.$$

*Then we have*

$$(8) \qquad \varphi'(0)\hat{\theta}_n - \sum_{i=1}^n \psi(e_i)\mathbf{z}_{i,n} = o_\mathbb{P}(1)$$

*and* $\hat{\theta}_n = O_\mathbb{P}(1)$. *Additionally, if the limit*

$$(9) \qquad \lim_{n \to \infty} \sum_{i=1}^{n-|k|} \mathbf{z}_{i,n} \mathbf{z}_{i+k,n}' = \Delta_k$$

*exists for each* $k \in \mathbb{Z}$, *then*

$$(10) \qquad \varphi'(0)\hat{\theta}_n \Rightarrow N(0, \Delta), \qquad \text{where } \Delta = \sum_{k \in \mathbb{Z}} \mathbb{E}[\psi(e_0)\psi(e_k)]\Delta_k.$$

Theorem 1 ensures the consistency of $\hat{\beta}_n$: $\hat{\theta}_n = O_\mathbb{P}(1)$ and $\Sigma_n^{-1} \to 0$ implies that $\hat{\beta}_n = o_\mathbb{P}(1)$. It is generally not trivial to establish the consistency of $M$-estimators. The convexity condition is quite useful in proving consistency; see [7, 23, 39] among others for regression models with independent errors. Recently, Berlinet, Liese and Vajda [9] considered consistency of $M$-estimates in regression models with strong mixing errors. This paper requires that the regressors $\mathbf{x}_i$ satisfy the condition that $n^{-1} \sum_{i=1}^n \delta_{\mathbf{x}_i}$ converges to some probability measure, where $\delta$ is the Dirac measure. This condition seems restrictive and it excludes some interesting cases (cf. Remark 1). For linear models with stationary causal errors, it is unclear how to establish the consistency and asymptotic normality without the convexity of $\rho$.

We now discuss condition (7). Since $\psi_i(\epsilon; \mathcal{F}_0) = \mathbb{E}[\psi(e_i + \epsilon)|\mathcal{F}_0]$ is the $i$th-step-ahead predicted mean, the quantity $\|\psi_i(\epsilon; \mathcal{F}_0) - \psi_i(\epsilon; \mathcal{F}_0^*)\| = \|\mathbb{E}[\psi(e_i + \epsilon)|\mathcal{F}_0] - \mathbb{E}[\psi(e_i^* + \epsilon)|\mathcal{F}_0^*]\|$ measures the contribution of $\varepsilon_0$ in predicting $\psi(e_i + \epsilon)$.

Hence (7) suggests short-range dependence in the sense that the cumulative contribution of $\varepsilon_0$ in predicting future values is finite. The following proposition provides a sufficient condition for (7). Recall that $F_i(\cdot|\mathcal{F}_0)$ is the conditional (or predictive) distribution function of $e_i$ given $\mathcal{F}_0$ and $f_i(\cdot|\mathcal{F}_0)$ is the conditional density. Let $\psi(u; \epsilon_0) = |\psi(u + \epsilon_0)| + |\psi(u - \epsilon_0)|$.

PROPOSITION 1. *Condition* (7) *holds under either* (i)

$$(11) \quad \sum_{i=1}^{\infty} \bar{\omega}(i) < \infty, \qquad \text{where } \bar{\omega}(i) = \int_{\mathbb{R}} \|f_i(u|\mathcal{F}_0) - f_i(u|\mathcal{F}_0^*)\| \psi(u; \epsilon_0) \, du,$$

*or* (ii) $\rho(x) = \rho_\alpha(x) = \alpha x^+ + (1-\alpha)(-x)^+$, $0 < \alpha < 1$, *and*

$$(12) \quad \sum_{i=1}^{\infty} \omega(i) < \infty, \qquad \text{where } \omega(i) = \sup_{u \in \mathbb{R}} \|F_i(u|\mathcal{F}_0) - F_i(u|\mathcal{F}_0^*)\|.$$

Proposition 1 easily follows from the identities $\mathbb{E}(\mathbf{1}_{e_i \leq u}|\mathcal{F}_0) = F_i(u|\mathcal{F}_0)$ and $\mathbb{E}[\psi(e_i + \epsilon)|\mathcal{F}_0] = \int_{\mathbb{R}} \psi(u + \epsilon) f_i(u|\mathcal{F}_0) \, du$. We omit the details.

In Proposition 1, case (ii) corresponds to quantile regression, an important non-least squares procedure. Condition (12) can be interpreted as follows. If the conditional distribution $[e_i|\mathcal{F}_0]$ does not depend on $\varepsilon_0$, then $F_i(u|\mathcal{F}_0) - F_i(u|\mathcal{F}_0^*) = 0$. The quantity $\sup_u \|F_i(u|\mathcal{F}_0) - F_i(u|\mathcal{F}_0^*)\|$ can thus be interpreted as the contribution of $\varepsilon_0$ in predicting $e_i$. In other words, $\sup_u \|F_i(u|\mathcal{F}_0) - F_i(u|\mathcal{F}_0^*)\|$ quantifies the degree of dependence of the predictive distribution $[e_i|\mathcal{F}_0]$ on $\varepsilon_0$. So (12) suggests that the cumulative contribution of $\varepsilon_0$ in predicting future values $(e_i)_{i \geq 1}$ is finite. Condition (11) delivers a similar message by incorporating the information of the target function $\psi = \rho'$ as weights into the distance between the two predictive distributions $[e_i|\mathcal{F}_0]$ and $[e_i^*|\mathcal{F}_0^*]$. To obtain Bahadur representations, stronger versions of (11) and (12) are needed; see (27) and (28).

REMARK 1. Many of the earlier results require that $\mathbf{x}_i$, $1 \leq i \leq n$, satisfy the condition that $\Sigma_n/n$ converges to a positive definite matrix ([8, 32] among others). This condition is not required in our setting. Consider the polynomial regression with design vectors $\mathbf{x}_i = (1, i, \ldots, i^{p-1})'$, $1 \leq i \leq n$. Then $\Sigma_n/n$ does not have a limit. Elementary but tedious calculations show that (A4) is satisfied and (9) holds with $\Delta_k = \text{Id}_p$. However, a condition of such type is needed in deriving strong Bahadur representations; see Theorem 3.

REMARK 2. In the expression of $\Delta$ in (10), the presence of the terms $\mathbb{E}[\psi(e_0)\psi(e_k)]\Delta_k$, $k \neq 0$, is due to the dependence of $(e_i)$.

REMARK 3. To apply Theorem 1 to quantile regression with $\rho(x) = \rho_\alpha(x)$, since $\psi(x) = \alpha - \mathbf{1}_{x \leq 0}$ and $\varphi(x) = \alpha - F(-x)$, we need to ensure that $e_i$ has a density at 0; see condition (A2). This problem is generally not easy. A simple sufficient condition is that the conditional density $f_1(\cdot|\mathcal{F}_0)$ exists. Without conditions of such type, the existence of marginal densities is not guaranteed. For example, consider the process $e_t = \sum_{i=0}^{\infty} 2\varepsilon_{t-i}/3^{i+1}$, where $\varepsilon_t$ are i.i.d. and $\mathbb{P}(\varepsilon_t = 1) = \mathbb{P}(\varepsilon_t = -1) = 1/2$. Then the conditional density does not exist and the marginal distribution does not have a density either. Solomyak [52] considered the absolute continuity of $\sum_{i=0}^{\infty} r^i \varepsilon_{t-i}$, $r \in (0,1)$.

2.2. *Bahadur representations.* Bahadur representations with appropriate rates are useful in the study of the asymptotic behavior of statistical estimators. For $M$-estimation under independent errors, various Bahadur representations have been derived; see, for example, [2, 4, 11, 26] and [47] among others. In particular, He and Shao [26] obtained a sharp almost sure bound under very general conditions on $\rho$. To obtain approximation rates for $M$-estimates of linear models with dependent errors, we need extra conditions on the behavior of the function $\psi_1(s; \mathcal{F}_0)$ at the neighborhood of $s = 0$:

(A5) There exists an $\epsilon_0 > 0$ such that

$$(13) \qquad L_i := \sup_{|s|,|t| \leq \epsilon_0, s \neq t} \frac{|\psi_1(s; \mathcal{F}_i) - \psi_1(t; \mathcal{F}_i)|}{|s-t|} \in \mathcal{L}^1.$$

(A6) Let $\psi_1(\cdot; \mathcal{F}_i) \in \mathcal{C}^l$, $l \geq 0$. For some $\epsilon_0 > 0$, $\sup_{|\epsilon| \leq \epsilon_0} \|\psi_1^{(l)}(\epsilon; \mathcal{F}_i)\| < \infty$ and

$$(14) \qquad \sum_{i=0}^{\infty} \sup_{|\epsilon| \leq \epsilon_0} \|\mathbb{E}[\psi_1^{(l)}(\epsilon; \mathcal{F}_i)|\mathcal{F}_0] - \mathbb{E}[\psi_1^{(l)}(\epsilon; \mathcal{F}_i^*)|\mathcal{F}_0^*]\| < \infty.$$

Condition (A5) suggests that the function $\psi_1(s; \mathcal{F}_i)$, $|s| \leq \epsilon_0$, is stochastically Lipschitz continuous at a neighborhood of 0. The function $\psi$ itself does not have to be Lipschitz continuous. Indeed, for $\psi(x) = \rho'_\alpha(x) = \alpha - \mathbf{1}_{x \leq 0}$, if the conditional density $f_1(\cdot|\mathcal{F}_0)$ is bounded, then (13) holds. For this $\psi$, we need to assume that the conditional density exists to ensure that $e_i$ has a density, which is a prerequisite for Bahadur representations for quantile estimates (cf. Remark 3). If $\psi \in \mathcal{C}^{l+1}$ satisfies $\sup_u |\psi^{(k)}(u)| < \infty$, $1 \leq k \leq l+1$, then (13) holds and a sufficient condition for (14) is that $\sum_{i=0}^{\infty} \|e_i - e_i^*\| < \infty$. In this case the existence of a conditional density is not required. Since $\mathbb{E}[\psi(e_1 - s)|\mathcal{F}_0] = \int_{\mathbb{R}} \psi(v) f_1(v + s|\mathcal{F}_0) \, dv$, by Fubini's theorem, a sufficient condition for (13) is $\int_{\mathbb{R}} |f'_1(u|\mathcal{F}_0)| \psi(u, \epsilon_0) \, du \in \mathcal{L}^1$. The latter holds if $\int_{\mathbb{R}} \|f'_1(u|\mathcal{F}_0)\| \psi(u, \epsilon_0) \, du < \infty$. The last condition (A6) is a generalization of (7). Section 2.3 gives sufficient conditions for (14).

Define $M$-processes $K_n(\theta) = \Omega_n(\theta) - \mathbb{E}[\Omega_n(\theta)]$ and $\tilde{K}(\beta) = \tilde{\Omega}_n(\beta) - \mathbb{E}[\tilde{\Omega}_n(\beta)]$, where

$$\Omega_n(\theta) = \sum_{i=1}^{n} \psi(e_i - \mathbf{z}'_{i,n}\theta)\mathbf{z}_{i,n} \quad \text{and} \quad \tilde{\Omega}_n(\beta) = \sum_{i=1}^{n} \psi(e_i - \mathbf{x}'_i\beta)\mathbf{x}_i,$$

(15)
$$\theta, \beta \in \mathbb{R}^p.$$

The $M$-process itself is an interesting subject of study and it plays an important role in $M$-estimation theory. Welsh [58] considered $M$-processes for linear models with i.i.d. errors. Theorems 2 and 3 present local oscillation rates for the $M$-processes $K_n$ and $\tilde{K}_n$. Corollary 1 provides a weak Bahadur representation for $\hat{\theta}_n$. Theorem 3 deals with $\tilde{K}$ and gives a strong Bahadur representation for $\hat{\beta}_n$.

THEOREM 2. *Assume* (A1)–(A5) *and assume* (A6) *holds with* $l = 0, \ldots, p$. *Let* $(\delta_n)_{n \in \mathbb{N}}$ *be a sequence of positive numbers such that*

(16)
$$\delta_n \to \infty \quad \text{and} \quad \delta_n r_n = \delta_n \max_{i \leq n} |\mathbf{z}_{i,n}| \to 0.$$

*Then*

(17)
$$\sup_{|\theta| \leq \delta_n} |K_n(\theta) - K_n(0)| = O_{\mathbb{P}}[\sqrt{\tau_n(\delta_n)} \log n + \delta_n \sqrt{\zeta_n(4)}],$$

*where*

(18)
$$\tau_n(\delta) = \sum_{i=1}^{n} |\mathbf{z}_{i,n}|^2 [m^2(|\mathbf{z}_{i,n}|\delta) + m^2(-|\mathbf{z}_{i,n}|\delta)], \qquad \delta > 0.$$

COROLLARY 1. *Assume* (A1)–(A5) *and assume* (A6) *holds with* $l = 0, \ldots, p$, *and* $\varphi(t) = t\varphi'(0) + O(t^2)$ *as* $t \to 0$. *Further assume* $\Omega(\hat{\theta}_n) = O_{\mathbb{P}}(r_n)$. *Then for any sequence* $c_n \to \infty$,

(19)
$$\varphi'(0)\hat{\theta}_n - \sum_{i=1}^{n} \psi(e_i)\mathbf{z}_{i,n} = O_{\mathbb{P}}[\sqrt{\tau_n(\delta_n)} \log n + \delta_n r_n],$$

*where* $\delta_n = \min(c_n, r_n^{-1/2})$.

*In particular, if as* $t \to 0$, $m(t) = O(|t|^\lambda)$ *for some* $\lambda > 0$, *then*

(20)
$$\varphi'(0)\hat{\theta}_n - \sum_{i=1}^{n} \psi(e_i)\mathbf{z}_{i,n} = O_{\mathbb{P}}[\sqrt{\zeta_n(2 + 2\lambda)} \log n + r_n].$$

REMARK 4. If $\psi$ is continuous, it is easily seen that the minimizer $\hat{\theta}_n$ solves the equation $\Omega(\hat{\theta}_n) = 0$. In the case that $\psi$ is discontinuous, the latter equation may not have a solution. To overcome this difficulty, in Corollary 1 we propose the approximate equation $\Omega(\hat{\theta}_n) = O_\mathbb{P}(r_n)$. An important example for discontinuous $\psi$ arises in quantile regression. Let $\psi(x) = \rho'_\alpha(x) = \alpha - \mathbf{1}_{x \leq 0}$, $0 < \alpha < 1$. The argument in Corollary 2 and Lemma 9 implies that the minimizer $\hat{\theta}_n$ satisfies $|\Omega(\hat{\theta}_n)| \leq (p+1)r_n$ almost surely.

THEOREM 3. (a) *Assume* (A1)–(A3), (A5) *and assume* (A6) *holds with* $l = 0, \ldots, p$. (b) *Let* $\lambda_n$ *be the minimum eigenvalue of* $\Sigma_n$. *Assume that*

$$\liminf_{n \to \infty} \lambda_n/n > 0, \qquad \xi_n(2) = O(n)$$

*and*

(21) $$\tilde{r}_n := \max_{j \leq n} |\mathbf{x}_j| = O[n^{1/2}(\log n)^{-2}].$$

*Let* $b_n = n^{-1/2}(\log n)^{3/2}(\log \log n)^{1/2+\iota}$, $\iota > 0$, $\bar{n} = 2^{\lceil \log n / \log 2 \rceil}$ *and* $q > 3/2$. *Then* (i)

(22) $$\sup_{|\beta| \leq b_n} |\tilde{K}_n(\beta) - \tilde{K}_n(0)| = O_{\text{a.s.}}(L_{\bar{n}} + B_{\bar{n}}),$$

*where* $B_n = b_n\sqrt{\xi_n(4)}(\log n)^{3/2}(\log \log n)^{(1+\iota)/2}$, $L_n = \sqrt{\tilde{\tau}_n(2b_n)}(\log n)^q$ *and*

(23) $$\tilde{\tau}_n(\delta) = \sum_{i=1}^{n} |\mathbf{x}_i|^2 [m^2(|\mathbf{x}_i|\delta) + m^2(-|\mathbf{x}_i|\delta)], \qquad \delta > 0.$$

*If additionally* $\varphi(t) = t\varphi'(0) + O(t^2)$ *and* $m(t) = O(\sqrt{t})$ *as* $t \to 0$ *and* $\tilde{\Omega}_n(\hat{\beta}_n) = O_{\text{a.s.}}(\tilde{r}_n)$, *then* (ii) $\hat{\beta}_n = O_{\text{a.s.}}(b_n)$ *and* (iii) *the strong Bahadur representation holds:*

(24) $$\varphi'(0)\Sigma_n\hat{\beta}_n - \sum_{i=1}^{n} \psi(e_i)\mathbf{x}_i = O_{\text{a.s.}}(L_{\bar{n}} + B_{\bar{n}} + \xi_n(3)b_n^2 + \tilde{r}_n).$$

COROLLARY 2. *Assume that* 0 *is the* $\alpha$*th quantile of* $e_i$, $f(0) > 0$ *and there is a constant* $C_0 < \infty$ *such that* $\sup_u f_1(u|\mathcal{F}_0) \leq C_0$, $f \in \mathcal{C}^1$, $\sup_{u \in \mathbb{R}} \|F_1^{(l)}(u|\mathcal{F}_i)\| < \infty$ *and*

(25) $$\sum_{i=0}^{\infty} \sup_{u \in \mathbb{R}} \|\mathbb{E}[F_1^{(l)}(u|\mathcal{F}_i)|\mathcal{F}_0] - \mathbb{E}[F_1^{(l)}(u|\mathcal{F}_i^*)|\mathcal{F}_0^*]\| < \infty, \qquad l = 0, \ldots, p.$$

*Assume that* $\mathbf{x}_i$ *satisfies conditions* (b) *in Theorem* 3. *Let* $\hat{\beta}_n$ *be a minimizer of* (2) *with* $\rho_\alpha$. *Then* (i) $\hat{\beta}_n = o_{\text{a.s.}}(b_n)$ *and* (ii) *the strong Bahadur representation holds:*

(26) $$f(0)\Sigma_n^{-1}\hat{\beta}_n - \sum_{i=1}^{n}(\alpha - \mathbf{1}_{e_i \leq 0})\mathbf{x}_i = O_{\text{a.s.}}(B_{\bar{n}} + L_{\bar{n}} + \xi_n(3)b_n^2 + \tilde{r}_n),$$

*where* $B_n = [b_n\xi_n^{1/2}(4)](\log n)^{3/2}(\log \log n)^{3/4}$ *and* $L_n = [b_n\xi_n(3)]^{1/2}(\log n)^q$.

We now discuss the bound in (24). Clearly the condition $\xi_n(2) = O(n)$ implies that $\tilde{r}_n = \max_{j \leq n} |\mathbf{x}_j| = O(n^{1/2})$. The condition on $\tilde{r}_n$ in (21) is not the weakest possible. For presentational clarity we adopt (21) since otherwise it involves quite tedious manipulations. If $\psi$ is continuous, then $\tilde{\Omega}_n(\hat{\beta}_n) = 0$. If additionally $\xi_n(4) = O(n)$ and $m(t) = O(|t|^\lambda)$, $1/2 \leq \lambda \leq 1$, the bound in (24) is $O[n^{(1-\lambda)/2}(\log n)^{q'}]$, $q' > 3$. For the bound in (26), elementary calculations show that, if $\xi_n(\kappa) = O(n)$, $2 < \kappa < 4$, then the bound is $O_{\text{a.s.}}(n^{1/\kappa})$; if $\xi_n(4) = O(n)$, then the bound becomes $O[n^{1/4}(\log n)^{q'}]$, $q' > 9/4$. The latter bound is optimal up to a multiplicative logarithmic factor since the classical [6] representation has the bound $O_{\text{a.s.}}[n^{-3/4}(\log \log n)^{3/4}]$.

If $\Sigma_n/n$ converges to a positive definite matrix $Q$ (say), then $\xi_n(2) = O(n)$ and the limit of $\lambda_n/n$ is the smallest eigenvalue of $Q$, which is strictly positive. To apply Theorems 2 and 3 and Corollary 1, we also need to verify $\varphi(t) = t\varphi'(0) + O(t^2)$ and know the order of magnitude of $m(\cdot)$; see the definitions of $\tau_n(\cdot)$ and $\tilde{\tau}_n(\cdot)$ by (18) and (23). Examples 1 and 2 below concern some commonly used $\rho$. Recall that $f$ is the density of $e_i$.

EXAMPLE 1. Assume (16). If $\psi$ has a derivative $\psi'$ satisfying $\sup_{|u| \leq \delta} \|\psi'(e_1 + u)\| < \infty$ for some $\delta > 0$, then $m(t) = O(|t|)$ as $t \to 0$ and $\tau_n(\delta_n) = O[\zeta_n(4)\delta_n^2]$. The latter claim easily follows from $\psi(e_1+t) - \psi(e_1) = \int_0^t \psi'(e_1+u)\,du$ and $m^2(t) \leq t \int_0^t \|\psi'(e_1+u)\|^2\,du = O(t^2)$. An important example is Huber's function $\rho(x) = (x^2 \mathbf{1}_{|x| \leq c})/2 + (c|x| - c^2/2)\mathbf{1}_{|x|>c}$, $c > 0$. Then $\psi(x) = \max[\min(x,c),-c]$ and $\varphi'(t) = \mathbb{P}(|e_1 + t| \leq c) = F(c-t) - F(-c-t)$. If $\sup_x f(x) < \infty$, then $\varphi(t) = t\varphi'(0) + O(t^2)$ and $m(t) = O(|t|)$ as $t \to 0$.

EXAMPLE 2 ($\mathcal{L}^q$-regression estimates). Assume (16). Let $\rho(t) = |t|^q$, $1 < q \leq 2$, and assume $\sup_v f(v) < \infty$. If $q \neq 3/2$, then $m(t) = O(|t|^{q'/2})$ and $\tau_n(\delta_n) = O[\zeta_n(2+q')\delta_n^{q'}]$, where $q' = \min(2, 2q-1)$. If $q = 3/2$, then $m(t) = O[|t|\log(1/|t|)]$ and $\tau_n(\delta_n) = \sum_{i=1}^n |\mathbf{z}_{i,n}|^4 (\log |\mathbf{z}_{i,n}|)^2 O(\delta_n^2)$. Here the bound of $m(t)$ follows from [2]. The bound of $\tau_n(\delta_n)$ when $q \neq 3/2$ can easily be obtained. If $q = 3/2$, since $r_n \delta_n \to 0$, then $|\log |\mathbf{z}_{i,n} \delta_n|| \leq 2|\log |\mathbf{z}_{i,n}||$ for sufficiently large $n$ and the stated bound for $\tau_n(\delta_n)$ follows.

If $\sup_x[f(x) + |f'(x)|] < \infty$ and $e_0 \in \mathcal{L}^{q-1}$, then $\varphi(t) = t\varphi'(0) + O(t^2)$. To this end note that $\psi(x) = q|x|^{q-1}\text{sgn}(x)$ and $\psi'(x) = q(q-1)|x|^{q-2}\text{sgn}(x)$. Let $|\delta| \leq 1$ and $D = \psi'(e_1 + \delta) - \psi'(e_1)$. If $|e_1| > 3$, then $|D| \leq |\delta|$. On the other hand,

$$\mathbb{E}(D\mathbf{1}_{|e_1| \leq 3}) = \int_{-3}^{3} \psi'(u)[f(u) - f(u-\delta)]\,du + \left[\int_{-3}^{-3+\delta} - \int_{3}^{3+\delta}\right]\psi'(u)f(u-\delta)\,du,$$

which is also of the order $O(\delta)$ since $\sup_x |f'(x)| < \infty$ and $\int_{-3}^{3} |\psi'(u)|\, du < \infty$. Therefore $\mathbb{E}[\int_0^t \psi'(e_1 + \delta) - \psi'(e_1)\, d\delta] = O(t^2)$, which implies $\varphi(t) - \varphi(0) = t\varphi'(0) + O(t^2)$.

2.3. *Sufficient conditions for* (14). Recall the projections $\mathcal{P}_k \cdot = \mathbb{E}(\cdot|\mathcal{F}_k) - \mathbb{E}(\cdot|\mathcal{F}_{k-1})$ and $\mathcal{F}_k^* = (\ldots, \varepsilon_{-1}, \varepsilon_0', \varepsilon_1, \ldots, \varepsilon_k)$. Proposition 2 provides sufficient conditions for (14) and (25). These sufficient conditions appear easy to work with; see applications in Section 3. Lemma 1 follows from Theorem 1 in [63].

LEMMA 1. *Assume that the process $X_t = g(\mathcal{F}_t) \in \mathcal{L}^2$. Let $g_n(\mathcal{F}_0) = \mathbb{E}[g(\mathcal{F}_n)|\mathcal{F}_0]$, $n \geq 0$. Then $\|\mathcal{P}_0 X_n\| \leq \|g(\mathcal{F}_n) - g(\mathcal{F}_n^*)\|$ and $\|\mathcal{P}_0 X_n\| \leq \|g_n(\mathcal{F}_0) - g_n(\mathcal{F}_0^*)\| \leq 2\|\mathcal{P}_0 X_n\|$.*

PROPOSITION 2. (i) *Assume that $f_1(\cdot|\mathcal{F}_i) \in \mathcal{C}^l$, $l \geq 0$, and*

$$\sum_{i=0}^{\infty} \bar{\omega}_l(i) < \infty,$$
(27)
$$\text{where } \bar{\omega}_l(i) = \int_{\mathbb{R}} \|f_1^{(l)}(u|\mathcal{F}_i) - f_1^{(l)}(u|\mathcal{F}_i^*)\|\psi(u;\epsilon_0)\, du.$$

*Then $\sum_{i=0}^{\infty} \sup_{|\epsilon|\leq \epsilon_0} \|\psi_1^{(l)}(\epsilon;\mathcal{F}_i) - \psi_1^{(l)}(\epsilon;\mathcal{F}_i^*)\| < \infty$ and (14) holds.*

(ii) *Let $\rho(x) = \rho_\alpha(x) = \alpha x^+ + (1-\alpha)(-x)^+$, $0 < \alpha < 1$. Then (25) holds if for $0 \leq l \leq p$*

$$\sum_{i=0}^{\infty} \omega_l(i) < \infty, \qquad \text{where } \omega_l(i) = \sup_{u\in\mathbb{R}} \|F_1^{(l)}(u|\mathcal{F}_i) - F_1^{(l)}(u|\mathcal{F}_i^*)\|.$$
(28)

PROOF. (i) Since $\psi_1(t;\mathcal{F}_i) = \int_\mathbb{R} \psi(v) f_1(v - t|\mathcal{F}_i)\, dv$, by Lemma 1, for $|t| \leq \epsilon_0$,

$$\|\mathbb{E}[\psi_1^{(l)}(t;\mathcal{F}_i)|\mathcal{F}_0] - \mathbb{E}[\psi_1^{(l)}(t;\mathcal{F}_i^*)|\mathcal{F}_0^*]\|$$
$$\leq 2\|\psi_1^{(l)}(t;\mathcal{F}_i) - \psi_1^{(l)}(t;\mathcal{F}_i^*)\|$$
$$= 2\left\|\int_\mathbb{R} \psi(v)[f_1^{(l)}(v - t|\mathcal{F}_i) - f_1^{(l)}(v - t|\mathcal{F}_i^*)]\, dv\right\|$$
$$\leq 2\int_\mathbb{R} |\psi(v)|\|f_1^{(l)}(v - t|\mathcal{F}_i) - f_1^{(l)}(v - t|\mathcal{F}_i^*)\|\, dv \leq 2\bar{\omega}_l(i).$$

So (i) holds.

(ii) This easily follows from Lemma 1. □

**3. Applications.** This section contains applications of results in Section 2 to linear models with errors being linear processes and some widely used nonlinear time series. For such processes the SRD conditions (27) and (28) can be verified.

3.1. *Linear processes.* Let $\varepsilon_i$ be i.i.d. random variables and $a_i$ real numbers such that

$$e_i = \sum_{j=0}^{\infty} a_j \varepsilon_{i-j} \tag{29}$$

exists. Without loss of generality let $a_0 = 1$. Let $F_\varepsilon$ be the distribution function of $\varepsilon_0$ and let $f_\varepsilon$ be its density. Propositions 3 and 4 provide simple sufficient conditions for (28) and (27), respectively. For $\gamma \in \mathbb{R}$ let the weighted measure $w_\gamma(du) = (1 + |u|)^\gamma \, du$. The proof of Proposition 4 is given in [64].

PROPOSITION 3. *Assume that $\varepsilon_0 \in \mathcal{L}^q$, $q > 0$, and that for some $C_0 < \infty$,*

$$\sup_u |f_\varepsilon^{(l)}(u)| < C_0, \qquad l = 0, \ldots, p. \tag{30}$$

*Then $\omega_l(i) = O(|a_i|^{q'/2})$, $q' = \min(2, q)$. Consequently (28) holds if*

$$\sum_{j=0}^{\infty} |a_j|^{q'/2} < \infty. \tag{31}$$

PROOF. Let $Z_n = \sum_{j=1}^{\infty} a_j \varepsilon_{n-j}$ and $Z_n^* = Z_n - a_n \varepsilon_0 + a_n \varepsilon_0'$. Then $F_1(u | \mathcal{F}_{n-1}) = F_\varepsilon(u - Z_n)$. By (30), since $\min(1, |x|^2) \le |x|^\delta$, $0 \le \delta \le 2$,

$$\begin{aligned}
\omega_l(n) &= \sup_u \|F_\varepsilon^{(l)}(u - Z_n) - F_\varepsilon^{(l)}(u - Z_n^*)\| \\
&\le \|\min(2C_0, C_0 |a_n \varepsilon_0 - a_n \varepsilon_0'|)\| \\
&\le 2C_0 \||a_n \varepsilon_0 - a_n \varepsilon_0'|^{q'/2}\| = O(|a_n|^{q'/2}).
\end{aligned} \tag{32}$$

The second assertion is obvious. □

PROPOSITION 4. *Let $1 < \gamma \le q$. Assume $\mathbb{E}(\varepsilon_0) = 0$, $\varepsilon_0 \in \mathcal{L}^q$, $\kappa_\gamma = \int_\mathbb{R} \psi^2(u) \times w_{-\gamma}(du) < \infty$ and*

$$\sum_{k=0}^{p+1} \int_\mathbb{R} |f_\varepsilon^{(k)}(v)|^2 w_\gamma(dv) < \infty. \tag{33}$$

*Then $\bar{\omega}_l(i) = O(|a_i|^{q'/2})$, $0 \le l \le p$, where $q' = \min(2, q)$, and (27) holds under (31).*

The condition $\kappa_\gamma < \infty$ controls the tails of $\psi$. Both Propositions 3 and 4 allow dependent and heavy-tailed errors. For the linear model (1) with heavy-tailed errors, it is more desirable to apply the $M$-estimation technique to estimate the unknown $\beta$ since the least squares procedure may result in estimators with erratic behavior. A popular model for such heavy-tailed processes is the moving average process $e_i = \sum_{j=0}^{\infty} a_j \varepsilon_{i-j}$, where $\varepsilon_i$ are i.i.d. random variables with stable distributions and $a_i$ are coefficients such that $e_i$ is well defined. Recently there has been a substantial interest in linear processes with heavy-tailed innovations; see [28, 54] and [61] among others. Davis, Knight and Liu [14] studied the behavior of the $M$-estimator in causal autoregressive models, while Davis and Wu [15] considered $M$-estimation in linear models. In the latter two papers the errors are assumed to be heavy-tailed, however, independent.

EXAMPLE 3. Let $\tau \geq 0$ be an integer and $\iota \in (1,2)$. Assume that the density $f_\varepsilon \in \mathcal{C}^\tau$ satisfies $f_\varepsilon^{(\tau)}(t) \sim |t|^{-1-\tau-\iota} h(|t|)$ as $t \to \pm\infty$, where $h$ is a slowly varying function [31], namely $\lim_{x \to \infty} h(x\lambda)/h(x) = 1$ for all $\lambda > 0$. By Karamata's theorem [31], simple calculations show that there exist real constants $C_j$, $0 \leq j < \tau$, such that $f_\varepsilon^{(j)}(t) \sim C_j |t|^{-j-1-\iota} h(|t|)$ as $t \to \infty$, and, for some constant $C$, $1 - F_\varepsilon(t) \sim C|t|^{-\iota} h(|t|)$ and $F_\varepsilon(-t) \sim C|t|^{-\iota} h(|t|)$ as $t \to \infty$. So $\varepsilon_i$ is in the domain of attraction of a stable distribution with index $\iota$ [31]. Let $\gamma < 1 + 2\iota$. Then $\int_\mathbb{R} [f_\varepsilon^{(j)}(t)]^2 w_\gamma(dt) < \infty$, $0 \leq j \leq \tau$. Hence (33) holds. As an interesting special case, let $\varepsilon_i$ be i.i.d. standard symmetric-$\alpha$-stable (S$\alpha$S) random variables with index $\iota \in (1,2)$. By Theorem 2.4.2 in [31], the density $f_\varepsilon(t) \sim c_\iota |t|^{-1-\iota}$ as $|t| \to \infty$, where $c_\iota = \pi^{-1} \zeta(1+\iota) \sin(\iota\pi/2)$. A similar argument shows that $f_\varepsilon^{(\tau)}(t) \sim c_{\iota,\tau} |t|^{-1-\tau-\iota}$, where $c_{\iota,\tau} = c_\iota \prod_{i=1}^\tau (-i - \iota)$.

EXAMPLE 4. If $e_i$ have finite variance, then (31) with $q'/2 = 1$ implies that the covariances are absolutely summable. It seems that the robust estimation problem of (1) with SRD linear process errors has been rarely studied in the literature. If (31) is barely violated, for example, if $a_n = n^{-\mu}$, $1/2 < \mu < 1$, then the errors are long-range dependent (LRD). In the LRD case, the $M$-estimates behave very differently from those in the i.i.d. or weakly dependent error cases in that they are asymptotically first-order equivalent, in probability, to the least squares estimate [37].

REMARK 5. The condition (31) seems almost necessary for the asymptotic normality of $\hat\theta_n$. Let $\varepsilon_i$ be i.i.d. standard S$\alpha$S random variables with index $\iota \in (1,2)$ and $a_n \sim n^{-\mu}$, $n \in \mathbb{N}$. Then (31) is reduced to $\iota\mu > 2$. Surgailis [54] showed that, if $\iota\mu < 2$, then the empirical process of $e_i$ satisfies a non-central limit theorem and the normalizing sequence is no longer $\sqrt{n}$. The asymptotic distribution of our estimate $\hat\theta_n$ is unknown when $\iota\mu < 2$.

3.2. *Nonlinear time series.* Many nonlinear time series models have the form

$$e_i = R(e_{i-1}, \varepsilon_i), \tag{34}$$

where $R$ is a measurable function and $\varepsilon_i$ are i.i.d. innovations. Diaconis and Freedman [16] showed that (34) has a stationary solution if for some $q > 0$ and $t_0$,

$$\mathbb{E}(\log \ell_{\varepsilon_0}) < 0 \quad \text{and} \quad \ell_{\varepsilon_0} + |R(t_0, \varepsilon_0)| \in \mathcal{L}^q, \tag{35}$$

$$\text{where } \ell_{\varepsilon_0} = \sup_{x \neq x'} \frac{|R(x, \varepsilon_0) - R(x', \varepsilon_0)|}{|x - x'|}.$$

In this case iterations of (34) lead to (3). Due to the Markovian structure of $(e_i)$, we can let $F_1(u|e_i) = F_1(u|\mathcal{F}_i)$ and $f_1(u|e_i) = f_1(u|\mathcal{F}_i)$ be the conditional distribution and density functions. Then $F_1(u|v) = \mathbb{P}[R(v, \varepsilon_i) \leq u]$ and $f_1(u|v) = \partial F_1(u|v)/\partial u$.

PROPOSITION 5. *Assume that there exists a constant $C_0 < \infty$ such that*

$$\sup_{u,v} \left| \frac{\partial F_1^{(l)}(u|v)}{\partial v} \right| + \sup_{u,v} |F_1^{(l)}(u|v)| < C_0, \qquad l = 0, \ldots, p. \tag{36}$$

*Then under* (35) *we have $\omega_l(i) = O(\chi^i)$ for some $\chi \in (0, 1)$ and hence* (28) *holds.*

PROOF. Let $(\varepsilon_i')_{i \in \mathbb{Z}}$ be an i.i.d. copy of $(\varepsilon_i)_{i \in \mathbb{Z}}$ and, for $i \geq 0$, $e_i^* = G(\ldots, \varepsilon_{-1}, \varepsilon_0', \varepsilon_1, \ldots, \varepsilon_i)$ and $e_i' = G(\ldots, \varepsilon_{-1}', \varepsilon_0', \varepsilon_1, \ldots, \varepsilon_i)$. By Theorem 2 in [65], under (35) there exist $\varsigma > 0$ and $\varrho \in (0, 1)$ such that $\|e_i' - e_i\|_\varsigma = O(\varrho^i)$. So $\|e_i^* - e_i\|_\varsigma = O(\|e_i' - e_i\|_\varsigma + \|e_i' - e_i^*\|_\varsigma) = O(\varrho^i)$. Observe that $F_1^{(l)}(u|\mathcal{F}_{n-1}) = F_1^{(l)}(u|e_{n-1})$. As in (32), by (36),

$$\|F_1^{(l)}(u|\mathcal{F}_{n-1}) - F_1^{(l)}(u|\mathcal{F}_{n-1}^*)\| = \|F_1^{(l)}(u|e_{n-1}) - F_1^{(l)}(u|e_{n-1}^*)\|$$

$$\leq \|\min(2C_0, C_0|e_{n-1} - e_{n-1}^*|)\|$$

$$\leq 2C_0 \||e_{n-1} - e_{n-1}^*|^{\min(1,\varsigma/2)}\| = O(\chi^n),$$

where $\chi = \varrho^{\min(1,\varsigma/2)}$. $\square$

EXAMPLE 5. Consider the autoregressive conditional heteroscedasticity (ARCH) model

$$e_i = \varepsilon_{i-1}\sqrt{a^2 + b^2 e_{i-1}^2},$$

where $\varepsilon_i$ are i.i.d. innovations and $a, b$ are real parameters such that

$$\mathbb{E}(\log|b\varepsilon_0|) < 0 \quad \text{and} \quad \varepsilon_0 \in \mathcal{L}^q \tag{37}$$

for some $q > 0$. Then $\ell_{\varepsilon_0} = |b\varepsilon_0|$ and (35) holds. Note that (37) imposes very mild moment conditions and it even allows $|e_i|$ to have infinite mean. Let $F_\varepsilon$ (resp. $f_\varepsilon$) be the distribution (resp. density) function of $\varepsilon_0$. Assume that $f_\varepsilon$ satisfies (30). Since $F_1(u|v) = F_\varepsilon(u/\sqrt{a^2 + b^2 v^2})$, simple calculations show that (30) implies (36). As an interesting special case, let $\varepsilon_i$ have the standard Student $t$-distribution with degrees of freedom $k > 0$. Then the density $f_\varepsilon(t) = [k/(k+t^2)]^{(1+k)/2}/C_k$, where $C_k = k^{1/2}B(1/2, k/2)$ and $B(\cdot, \cdot)$ is the beta function. Clearly (30) and (35) are satisfied. In certain applications it is desirable to use ARCH models with Student-$t$ innovations to allow heavy tails [56].

Proposition 6 below gives a sufficient condition for (27) for the process

$$(38) \qquad e_i = \nu(e_{i-1}) + \varepsilon_i,$$

where $\nu$ is a Lipschitz continuous function such that the Lipschitz constant

$$(39) \qquad \ell_\nu := \sup_{a \neq b} \frac{|\nu(a) - \nu(b)|}{|a - b|} < 1$$

and $\mathbb{E}(|\varepsilon_i|^\alpha) < \infty$ for some $\alpha > 1$. The condition $\ell_\nu < 1$ implies that the nonlinear time series (38) has a unique stationary distribution. A prominent example of (38) is the threshold autoregressive process $e_{i+1} = \alpha_1 e_i^+ + \alpha_2(-e_i)^+ + \varepsilon_{i+1}$, where $\alpha_1, \alpha_2$ are real coefficients [55]. In this example (39) is satisfied if $\max(|\alpha_1|, |\alpha_2|) < 1$. If the process (34) is of the form (38), then condition (27) can be simplified.

PROPOSITION 6. *Assume $\int_\mathbb{R} \psi^2(t) w_{-\gamma}(dt) < \infty$ for some $\gamma > 1$. Further assume (39), $\varepsilon_i \in \mathcal{L}^q$, $\gamma < q < \gamma + 2$ and that $f_\varepsilon$ satisfies (33). Then there exists $\chi \in (0, 1)$ such that*

$$(40) \qquad \bar{\omega}_l(i) = O(\chi^i), \qquad 0 \leq l \leq p.$$

REMARK 6. In Proposition 6, $\varepsilon_i$ are allowed to have infinite variances. Proposition 6 follows from Theorem 2 in [65]. For a proof see [64].

## 4. Proofs of results in Section 2.2.

LEMMA 2. *Let $T_n = \sum_{i=1}^n \psi(e_i) \mathbf{z}_{i,n}$.*

(i) *Assume $\mathbb{E}[\psi(e_i)] = 0$, $\|\psi(e_i)\| < \infty$ and*

$$(41) \qquad \sum_{i=1}^\infty \|\mathbb{E}[\psi(e_i)|\mathcal{F}_0] - \mathbb{E}[\psi(e_i^*)|\mathcal{F}_0^*]\| < \infty.$$

*Then $\|T_n\| = O(1)$.*

(ii) *If in addition* (9) *and* (A4) *hold, then* $T_n \Rightarrow N(0, \Delta)$.

PROOF. (i) For $k \geq 0$ let $J_k = \sum_{i=1}^{n} \mathcal{P}_{i-k}\psi(e_i)\mathbf{z}_{i,n}$. Note that the summands of $J_k$ are martingale differences. By the orthogonality, since $\sum_{i=1}^{n} \mathbf{z}_{i,n}\mathbf{z}'_{i,n} = \mathrm{Id}_p$ and $k \geq 0$,

$$(42) \quad \|J_k\|^2 = \sum_{i=1}^{n} \|\mathcal{P}_{i-k}\psi(e_i)\mathbf{z}_{i,n}\|^2 = \sum_{i=1}^{n} |\mathbf{z}_{i,n}|^2 \|\mathcal{P}_0\psi(e_k)\|^2 = p\|\mathcal{P}_0\psi(e_k)\|^2.$$

By Lemma 1, $\|\mathcal{P}_0\psi(e_k)\| \leq \|\mathbb{E}[\psi(e_k)|\mathcal{F}_0] - \mathbb{E}[\psi(e_k^*)|\mathcal{F}_0^*]\|$. By (41), $\sum_{k=0}^{\infty} \|J_k\| < \infty$ and consequently $\|T_n\| = O(1)$ since $T_n = \sum_{k=0}^{\infty} J_k$.

(ii) We now show $T_n \Rightarrow N(0, \Delta)$. Let $\mathbf{c}$ be a $p$-dimensional column vector with $|\mathbf{c}| = 1$ and $u_{i,n} = \mathbf{c}'\mathbf{z}_{i,n}$. By the Cramér–Wold device, it suffices to verify that

$$(43) \quad \sum_{i=1}^{n} u_{i,n}\psi(e_i) \Rightarrow N(0, \mathbf{c}'\Delta\mathbf{c}).$$

Since $\sum_{i=1}^{n} \mathbf{z}_{i,n}\mathbf{z}'_{i,n} = \mathrm{Id}_p$, $d_n := (\sum_{i=1}^{n} u_{i,n}^2)^{1/2} = 1$. By (A4),

$$(44) \quad \lim_{n \to \infty} \frac{\max_{i \leq n} |u_{i,n}|}{d_n} = 0.$$

By (9), for each $k \geq 0$,

$$(45) \quad \lim_{n \to \infty} \frac{\sum_{i=1}^{n-k} u_{i,n} u_{i+k,n}}{d_n^2} = \mathbf{c}'\Delta_k\mathbf{c}.$$

Write $\psi(e_i) = \sum_{j=0}^{\infty} \alpha_j \eta_{i,i-j}$, where $\alpha_j = \|\mathcal{P}_{i-j}\psi(e_i)\|$ and $\eta_{i,i-j} = \mathcal{P}_{i-j}\psi(e_i)/\alpha_j$. Then (41) entails $\sum_{j=0}^{\infty} \alpha_j < \infty$. By the argument in the proof of Theorem 1(i) in [25], (44) and (45) imply (43). (Theorem 1(i) in [25] is not yet directly applicable: condition (5,a) therein requires $d_n \to \infty$ and condition (7) therein requires $\sum_{k \in \mathbb{Z}} \mathbf{c}'\Delta_k\mathbf{c} > 0$. However, a careful examination of Hannan's [25] proof indicates that his conditions (5,a) and (7) are not needed in deriving (43) from (44) and (45). Also note that there is a typo in [25], (5,b). The correct version should be of the form (44).) □

LEMMA 3. *Let* $(\alpha_{ni})_{i=1}^{n}$, $n \in \mathbb{N}$, *be a triangular array of real numbers such that* $\sum_{i=1}^{n} \alpha_{ni}^2 \leq 1$ *and* $\varpi_n := \max_{i \leq n} |\alpha_{ni}| \to 0$. *Assume* (A1), (A3) *and* (7). *Let* $\eta_{n,i} = \rho(e_i - \alpha_{ni}) - \rho(e_i) + \alpha_{ni}\psi(e_i)$. *Then* $\mathrm{var}(\sum_{i=1}^{n} \eta_{n,i}) \to 0$.

PROOF. Let $I \in \mathbb{N}$ and $0 < \mu < \epsilon_0$. By (A3), (7) and Lemma 1, we have

$$(46) \quad \sum_{k=0}^{\infty} \sup_{|\epsilon| \leq \mu} \|\mathcal{P}_0[\psi(e_k) - \psi(e_k - \epsilon)]\| \leq \sum_{k=0}^{I} \sup_{|\epsilon| \leq \mu} \|\psi(e_k) - \psi(e_k - \epsilon)\|$$
$$+ \sum_{k=1+I}^{\infty} 2 \sup_{|\epsilon| \leq \mu} \|\mathcal{P}_0\psi(e_k - \epsilon)\| \to 0$$

by first letting $\mu \to 0$ and then $I \to \infty$. Let $A_n = \sum_{i=1}^{n} \alpha_{ni}^2$ and

$$Z_{k,n}(t) = \sum_{i=1}^{n} \mathcal{P}_{i-k}[\psi(e_i) - \psi(e_i - \alpha_{ni}t)]\alpha_{ni}, \qquad 0 \le t \le 1, \ k \ge 0.$$

Note that the summands of $Z_k(t)$ are martingale differences. Since $\varpi_n \to 0$,

(47)
$$\sup_{0 \le t \le 1} \|Z_{k,n}(t)\|^2 = \sup_{0 \le t \le 1} \sum_{i=1}^{n} \alpha_{ni}^2 \|\mathcal{P}_{i-k}[\psi(e_i) - \psi(e_i - \alpha_{ni}t)]\|^2$$
$$\le A_n \sup_{|\epsilon| \le \mu} \|\mathcal{P}_0[\psi(e_k) - \psi(e_k - \epsilon)]\|^2$$

holds for large $n$. Since $\eta_{n,i} = \int_0^1 [\psi(e_i) - \psi(e_i - t\alpha_{ni})]\alpha_{ni}\,dt$ and $A_n \le 1$,

$$\left\|\sum_{i=1}^{n}(\eta_{n,i} - \mathbb{E}\eta_{n,i})\right\| = \left\|\sum_{k=0}^{\infty}\int_0^1 Z_{k,n}(t)\,dt\right\| \le \sum_{k=0}^{\infty}\int_0^1 \|Z_{k,n}(t)\|\,dt,$$

which by (46) and (47) converges to 0 as $n \to \infty$ and $\mu \downarrow 0$. □

PROPOSITION 7. *Under conditions of Theorem* 1, *we have for any* $c > 0$ *that*

(48) $$D_n(c) := \sup_{|\theta| \le c}\left|\sum_{i=1}^{n}[\rho(e_i - \mathbf{z}_{i,n}'\theta) - \rho(e_i) + \mathbf{z}_{i,n}'\theta\psi(e_i)] - \frac{\varphi'(0)}{2}|\theta|^2\right| \to 0$$

*in probability.*

PROOF. We should use the argument in [7]. Let $\eta_i(\theta) = \rho(e_i - \mathbf{z}_{i,n}'\theta) - \rho(e_i) + \mathbf{z}_{i,n}'\theta\psi(e_i)$. For a fixed vector $\theta$ with $|\theta| \le c$, let $\alpha_{ni} = \mathbf{z}_{i,n}'\theta$. Then $\sum_{i=1}^{n}\alpha_{ni}^2 \le c^2$ and $\max_{i \le n}|\alpha_{ni}| \le cr_n \to 0$. By Lemma 3, $\text{var}[\sum_{i=1}^{n}\eta_i(\theta)] \to 0$. Note that $\sum_{i=1}^{n}\mathbf{z}_{i,n}\mathbf{z}_{i,n}' = \text{Id}_p$. By Lemma 1 in [7], under (A1) and (A2) the bias

$$\sum_{i=1}^{n}\mathbb{E}[\eta_i(\theta)] - \frac{\varphi'(0)}{2}|\theta|^2 = \sum_{i=1}^{n}\left[\frac{\varphi'(0)}{2}|\mathbf{z}_{i,n}'\theta|^2 + o(|\mathbf{z}_{i,n}'\theta|^2)\right] - \frac{\varphi'(0)}{2}|\theta|^2$$
$$= \frac{\varphi'(0)}{2}\sum_{i=1}^{n}[\theta'\mathbf{z}_{i,n}\mathbf{z}_{i,n}'\theta + o(|\mathbf{z}_{i,n}'\theta|^2)] - \frac{\varphi'(0)}{2}|\theta|^2$$
$$= o[\zeta_n(2)].$$

So $\sum_{i=1}^{n}\eta_i(\theta) \to \varphi'(0)|\theta|^2/2$ pointwise. Since $\eta_i(\theta)$, $1 \le i \le n$, are convex functions of $\theta$, by the convexity lemma in ([42], page 187), we have uniform in-probability convergence. (A nonstochastic version is given in [48], Theorem 10.8, page 90. A subsequencing argument leads to the in-probability-convergence version. See also Appendix II in [1] and [7] for more details.) □

PROOF OF THEOREM 1. The relation (8) easily follows from properties of convex functions; see, for example, the proofs of Theorems 2.2 and 2.4 in [7] and Theorem 1 in [42]. We omit the details. A proof is given in [64]. That $\hat{\theta}_n = O_{\mathbb{P}}(1)$ follows easily from Lemma 2 and (8). If (9) holds, again by Lemma 2 we have the central limit theorem (10). □

PROOF OF THEOREM 2. Write $K_n = M_n + N_n$, where

$$(49) \qquad M_n(\theta) = \sum_{i=1}^n \{\psi(e_i - \mathbf{z}'_{i,n}\theta) - \mathbb{E}[\psi(e_i - \mathbf{z}'_{i,n}\theta)|\mathcal{F}_{i-1}]\}\mathbf{z}_{i,n}$$

and, noting that $\mathbb{E}[\psi(e_i - \mathbf{z}'_{i,n}\theta)|\mathcal{F}_{i-1}] = \psi_1(-\mathbf{z}'_{i,n}\theta; \mathcal{F}_{i-1})$,

$$(50) \qquad N_n(\theta) = \sum_{i=1}^n \{\psi_1(-\mathbf{z}'_{i,n}\theta; \mathcal{F}_{i-1}) - \varphi(-\mathbf{z}'_{i,n}\theta)\}\mathbf{z}_{i,n}.$$

The summands of $M_n$ form (triangular array) martingale differences with respect to the filter $\sigma(\mathcal{F}_i)$. Since $\sum_{i=1}^n \mathbf{z}_{i,n}\mathbf{z}'_{i,n} = \mathrm{Id}_p$, we have $n^{-1/2} = O(r_n)$ and by (16), $n^{-2} = O(r_n^4) = O[\delta_n\zeta_n(4)]$. By Lemmas 4 and 5, (17) follows since $n^{-3} = O[\delta_n\zeta_n(4)]$. □

LEMMA 4. *Assume* (A5) *and* (16). *Then*

$$(51) \qquad \sup_{|\theta|\leq\delta_n} |M_n(\theta) - M_n(0)| = O_{\mathbb{P}}[\sqrt{\tau_n(\delta_n)}\log n + n^{-3}].$$

PROOF. Since $p = \sum_{i=1}^n \mathbf{z}'_{i,n}\mathbf{z}_{i,n} \leq nr_n^2$, (16) implies that $\delta_n = o(\sqrt{n})$. It suffices to show that the left-hand side of (51) has the bound $O_{\mathbb{P}}[g_n\sqrt{\tau_n(\delta_n)}\log n + n^{-3}]$ for any positive sequence $g_n \to \infty$. Assume that $g_n \geq 3$ for all $n$. Let

$$\phi_n = 2g_n\sqrt{\tau_n(\delta_n)}\log n, \qquad t_n = \frac{g_n\sqrt{\tau_n(\delta_n)}}{\log g_n}, \qquad u_n = t_n^2,$$

$$\eta_i(\theta) = [\psi(e_i - \mathbf{z}'_{i,n}\theta) - \psi(e_i)]\mathbf{z}_{i,n}, \qquad T_n = \max_{i\leq n}\sup_{|\theta|\leq\delta_n}|\eta_i(\theta)|,$$

$$U_n = \sum_{i=1}^n \mathbb{E}\{[\psi(e_i + |\mathbf{z}_{i,n}|\delta_n) - \psi(e_i - |\mathbf{z}_{i,n}|\delta_n)]^2|\mathcal{F}_{i-1}\}|\mathbf{z}_{i,n}|^2.$$

Since $\psi$ is monotone, for $\delta \geq 0$,

$$\sup_{|\theta|\leq\delta}|\eta_i(\theta)| \leq |\mathbf{z}_{i,n}|\max[|\psi(e_i - |\mathbf{z}_{i,n}|\delta) - \psi(e_i)|, |\psi(e_i + |\mathbf{z}_{i,n}|\delta) - \psi(e_i)|]$$

$$\leq |\mathbf{z}_{i,n}|[\psi(e_i + |\mathbf{z}_{i,n}|\delta) - \psi(e_i - |\mathbf{z}_{i,n}|\delta)].$$

So $\mathbb{E}[\sup_{|\theta|\leq \delta}|\eta_i(\theta)|^2] \leq 2|\mathbf{z}_{i,n}|^2[m^2(-|\mathbf{z}_{i,n}|\delta) + m^2(|\mathbf{z}_{i,n}|\delta)]$, $\mathbb{E}(T_n^2) \leq 2\tau_n(\delta_n)$ and

(52) $\quad \mathbb{P}(T_n \geq t_n) \leq t_n^{-2}\mathbb{E}(T_n^2) \leq 2t_n^{-2}\tau_n(\delta_n) = O[(g_n^{-1}\log g_n)^2] \to 0.$

Similarly, $\mathbb{E}(U_n) \leq 2\tau_n(\delta_n)$ and

(53) $\quad \mathbb{P}(U_n \geq u_n) \leq u_n^{-1}\mathbb{E}(U_n) = O[(g_n^{-1}\log g_n)^2] \to 0.$

Write $\mathbf{z}_{i,n} = (z_{i1,n}, \ldots, z_{ip,n})'$. For notational simplicity we write $z_{ij}$ for $z_{ij,n}$, $1 \leq j \leq p$. Let $\Pi_p = \{-1,+1\}^p$. For $i \in \mathbb{N}$ let $D_{\mathbf{z}}(i) = (2 \times \mathbf{1}_{z_{i1}\geq 0} - 1, \ldots, 2 \times \mathbf{1}_{z_{ip}\geq 0} - 1) \in \Pi_p$. For $\mathbf{d} \in \Pi_p$ and $1 \leq j \leq p$ define

$$M_{n,j,\mathbf{d}}(\theta) = \sum_{i=1}^n \{\psi(e_i - \mathbf{z}'_{i,n}\theta) - \mathbb{E}[\psi(e_i - \mathbf{z}'_{i,n}\theta)|\mathcal{F}_{i-1}]\}z_{ij}\mathbf{1}_{D_{\mathbf{z}}(i)=\mathbf{d}}.$$

Since $M_n = \sum_{\mathbf{d}\in\Pi_p}(M_{n,1,\mathbf{d}}, \ldots, M_{n,p,\mathbf{d}})'$, it suffices to show that (51) holds with $M_n$ therein replaced by $M_{n,j,\mathbf{d}}$ for all $\mathbf{d} \in \Pi_p$ and $1 \leq j \leq p$. To this end, for presentational clarity we consider $j = 1$ and $\mathbf{d} = (1, -1, 1, 1, \ldots, 1)$. The other cases similarly follow.

Let $|\theta| \leq \delta_n$, $\eta_{i,j,\mathbf{d}}(\theta) = [\psi(e_i - \mathbf{z}'_{i,n}\theta) - \psi(e_i)]z_{ij}\mathbf{1}_{D_{\mathbf{z}}(i)=\mathbf{d}}$ and

$$B_n(\theta) = \sum_{i=1}^n \mathbb{E}[\eta_{i,j,\mathbf{d}}(\theta)\mathbf{1}_{|\eta_{i,j,\mathbf{d}}(\theta)|>t_n}|\mathcal{F}_{i-1}].$$

Then for large $n$, since $u_n = o(t_n\phi_n)$,

(54) $\quad \mathbb{P}(|B_n(\theta)| \geq \phi_n, U_n \leq u_n) \leq \mathbb{P}\left[t_n^{-1}\sum_{i=1}^n \mathbb{E}[\eta_{i,j,\mathbf{d}}^2(\theta)|\mathcal{F}_{i-1}] \geq \phi_n, U_n \leq u_n\right]$
$\leq \mathbb{P}(t_n^{-1}U_n \geq \phi_n, U_n \leq u_n) = 0.$

Since $\mathcal{P}_i[\eta_{i,j,\mathbf{d}}(\theta)\mathbf{1}_{|\eta_{i,j,\mathbf{d}}(\theta)|\leq t_n}]$, $i = 1, \ldots, n$, form bounded martingale differences, by Proposition 2.1 in [19] and (54), for $|\theta| \leq \delta_n$,

$\mathbb{P}[|M_{n,j,\mathbf{d}}(\theta) - M_{n,j,\mathbf{d}}(0)| \geq 2\phi_n, T_n \leq t_n, U_n \leq u_n]$

$\leq \mathbb{P}\left[\left|\sum_{i=1}^n \mathcal{P}_i[\eta_{i,j,\mathbf{d}}(\theta)\mathbf{1}_{|\eta_{i,j,\mathbf{d}}(\theta)|\leq t_n}]\right| \geq \phi_n, T_n \leq t_n, U_n \leq u_n\right]$

(55) $\quad + \mathbb{P}\left[\left|\sum_{i=1}^n \mathcal{P}_i[\eta_{i,j,\mathbf{d}}(\theta)\mathbf{1}_{|\eta_{i,j,\mathbf{d}}(\theta)|>t_n}]\right| \geq \phi_n, T_n \leq t_n, U_n \leq u_n\right]$

$= O[\exp\{-\phi_n^2/(4t_n\phi_n + 2u_n)\}] + \mathbb{P}(|B_n(\theta)| \geq \phi_n, T_n \leq t_n, U_n \leq u_n)$

$= O[\exp\{-\phi_n^2/(4t_n\phi_n + 2u_n)\}].$

Let $\ell = n^8$ and $G_\ell = \{(k_1/\ell, \ldots, k_p/\ell) : k_i \in \mathbb{Z}, |k_i| \leq n^9\}$. Note that $G_\ell$ has $(2n^9 + 1)^p$ points. By (55), since $t_n \phi_n \log n = o(\phi_n^2)$ and $u_n \log n = o(\phi_n^2)$, for any $\varsigma > 1$ we have

$$
(56) \quad \begin{aligned}
\mathbb{P}\biggl[\sup_{\theta \in G_\ell} |M_{n,1,\mathbf{d}}(\theta) - M_{n,1,\mathbf{d}}(0)| &\geq 2\phi_n, T_n \leq t_n, U_n \leq u_n\biggr] \\
&= O(n^{9p})O[\exp\{-\phi_n^2/(4t_n\phi_n + 2u_n)\}] = O(n^{-\varsigma p}),
\end{aligned}
$$

which by (52) and (53) implies

$$
(57) \quad \lim_{n \to \infty} \mathbb{P}\biggl[\sup_{\theta \in G_\ell} |M_{n,1,\mathbf{d}}(\theta) - M_{n,1,\mathbf{d}}(0)| \geq 2\phi_n\biggr] = 0.
$$

For $a \in \mathbb{R}$ let $\lceil a \rceil_\ell = \lceil a\ell \rceil/\ell$ and $\lfloor a \rfloor_\ell = \lfloor a\ell \rfloor/\ell$. Write $\langle a \rangle_{\ell,1} = \lfloor a \rfloor_\ell$ and $\langle a \rangle_{\ell,-1} = \lceil a \rceil_\ell$. Let $\mathbf{d} = (d_1, \ldots, d_p) \in \Pi_p$. For a vector $\theta = (\theta_1, \ldots, \theta_p)'$ let $\langle \theta \rangle_{\ell,\mathbf{d}} = (\langle \theta_1 \rangle_{\ell,d_1}, \ldots, \langle \theta_p \rangle_{\ell,d_p})$. For example, for $\mathbf{d} = (1, -1, 1, 1, \ldots, 1)$, we have $\langle \theta \rangle_{\ell,\mathbf{d}} = (\lfloor \theta_1 \rfloor_\ell, \lceil \theta_2 \rceil_\ell, \lfloor \theta_3 \rfloor_\ell, \ldots, \lfloor \theta_p \rfloor_\ell)$ and $\langle \theta \rangle_{\ell,-\mathbf{d}} = (\lceil \theta_1 \rceil_\ell, \lfloor \theta_2 \rfloor_\ell, \lceil \theta_3 \rceil_\ell, \ldots, \lceil \theta_p \rceil_\ell)$. Observe that for this $\mathbf{d}$, $\eta_{i,1,\mathbf{d}}(\langle \theta \rangle_{\ell,-\mathbf{d}}) \leq \eta_{i,1,\mathbf{d}}(\theta) \leq \eta_{i,j,\mathbf{d}}(\langle \theta \rangle_{\ell,\mathbf{d}})$ since $\psi$ is nondecreasing.

Let $|s|, |t| \leq r_n \delta_n$. By (13), $|\mathbb{E}[\psi(e_i - t) - \psi(e_i - s)|\mathcal{F}_{i-1}]| \leq L_{i-1}|s-t|$ for all large $n$ since $r_n \delta_n \to 0$. Let $V_n = \sum_{i=1}^n L_{i-1}$. Again by (13), $\mathbb{P}(V_n \geq n^4) \leq n^{-4}\mathbb{E}(V_n) = O(n^{-3})$. Since $|\theta - \langle \theta \rangle_{\ell,\mathbf{d}}| = O(\ell^{-1})$, we have $\max_{i \leq n} |\mathbf{z}'_{i,n}(\theta - \langle \theta \rangle_{\ell,\mathbf{d}})| = o(\ell^{-1})$, and by (16),

$$
\sup_{|\theta| \leq \delta_n} \sum_{i=1}^n |\mathbb{E}[\eta_i(\langle \theta \rangle_{\ell,\mathbf{d}}) - \eta_i(\theta)|\mathcal{F}_{i-1}]| \leq C\ell^{-1}V_n.
$$

Therefore, for all $|\theta| \leq \delta_n$,

$$
(58) \quad \begin{aligned}
M_{n,1,\mathbf{d}}(\langle \theta \rangle_{\ell,-\mathbf{d}}) - M_{n,1,\mathbf{d}}(0) - CV_n/\ell \\
\leq M_{n,1,\mathbf{d}}(\theta) - M_{n,1,\mathbf{d}}(0) \\
\leq M_{n,1,\mathbf{d}}(\langle \theta \rangle_{\ell,\mathbf{d}}) - M_{n,1,\mathbf{d}}(0) + CV_n/\ell,
\end{aligned}
$$

which implies (51) in view of (57) and $V_n/\ell = o_\mathbb{P}(n^4/\ell) = o_\mathbb{P}(n^{-3})$.  □

LEMMA 5. *Assume that* $(\delta_n)$ *satisfies* (16) *and* (A6) *holds with* $l = 0, \ldots, p$. *Then*

$$
(59) \quad \biggl\| \sup_{|\mathbf{g}| \leq \delta_n} |N_n(\mathbf{g}) - N_n(0)| \biggr\| = O[\sqrt{\zeta_n(4)}\delta_n].
$$

PROOF. Let $I = \{\alpha_1, \ldots, \alpha_q\} \subseteq \{1, \ldots, p\}$ be a nonempty set and $1 \leq \alpha_1 < \cdots < \alpha_q$. For a $p$-dimensional vector $\mathbf{u} = (u_1, \ldots, u_p)$ let $\mathbf{u}_I = (u_1 \mathbf{1}_{1 \in I}, \ldots, u_p \mathbf{1}_{p \in I})$. Note that the $j$th component of $\mathbf{u}_I$ is 0 if $j \notin I$, $1 \leq j \leq p$. Write

$$
\int_0^{\mathbf{g}_I} \frac{\partial^q N_n(\mathbf{u}_I)}{\partial \mathbf{u}_I} d\mathbf{u}_I = \int_0^{g_{\alpha_1}} \cdots \int_0^{g_{\alpha_q}} \frac{\partial^q N_n(\mathbf{u}_I)}{\partial u_{\alpha_1} \cdots \partial u_{\alpha_q}} du_{\alpha_1} \cdots du_{\alpha_q}
$$

and $\mathbf{w}_i = \mathbf{z}_{i,n} z_{i,\alpha_1} \cdots z_{i,\alpha_q}$. Observe that

$$\left|\frac{\partial^q N_n(\mathbf{u}_I)}{\partial \mathbf{u}_I}\right| = \left|\sum_{i=1}^n [\psi_1^{(q)}(-\mathbf{z}'_{i,n}\mathbf{u}_I; \mathcal{F}_{i-1}) - \varphi^{(q)}(-\mathbf{z}'_{i,n}\mathbf{u}_I)]\mathbf{w}_i\right|.$$

Let $|\mathbf{u}| \leq p\delta_n$ and $k \in \mathbb{N}$. Since $\max_{i \leq n} |\mathbf{z}_{i,n}\mathbf{u}| \leq pr_n\delta_n \to 0$, by Lemma 1, for large $n$,

$$\left\|\sum_{i=1}^n \mathcal{P}_{i-k}\psi_1^{(q)}(-\mathbf{z}'_{i,n}\mathbf{u}_I; \mathcal{F}_{i-1})\mathbf{w}_i\right\|^2$$

$$= \sum_{i=1}^n |\mathbf{w}_i|^2 \|\mathcal{P}_{i-k}\psi_1^{(q)}(-\mathbf{z}'_{i,n}\mathbf{u}_I; \mathcal{F}_{i-1})\|^2$$

$$\leq \sum_{i=1}^n |\mathbf{w}_i|^2 \sup_{|\epsilon| \leq \epsilon_0} \|\mathbb{E}[\psi_1^{(q)}(\epsilon; \mathcal{F}_{k-1})|\mathcal{F}_0] - \mathbb{E}[\psi_1^{(q)}(\epsilon; \mathcal{F}_{k-1}^*)|\mathcal{F}_0^*]\|^2.$$

As in the proof of (i) of Lemma 2, if (A6) holds with $l = 0, \ldots, p$, then $\|\partial^q N_n(\mathbf{u}_I)/\partial \mathbf{u}_I\| = O[\zeta_n^{1/2}(2 + 2q)]$ uniformly over $|\mathbf{u}| \leq p\delta_n$. Consequently,

$$\left\|\sup_{|\mathbf{g}| \leq \delta_n} \int_0^{\mathbf{g}_I} \left|\frac{\partial^q N_n(\mathbf{u}_I)}{\partial \mathbf{u}_I}\right| d\mathbf{u}_I\right\| \leq \left\|\int_{-\delta_n}^{\delta_n} \cdots \int_{-\delta_n}^{\delta_n} \left|\frac{\partial^q N_n(\mathbf{u}_I)}{\partial \mathbf{u}_I}\right| d\mathbf{u}_I\right\|$$

(60)
$$\leq \int_{-\delta_n}^{\delta_n} \cdots \int_{-\delta_n}^{\delta_n} \left\|\frac{\partial^q N_n(\mathbf{u}_I)}{\partial \mathbf{u}_I}\right\| d\mathbf{u}_I$$

$$= O[\delta_n^q \zeta_n^{1/2}(2 + 2q)].$$

By (16), $\delta_n^q \sqrt{\zeta_n(2+2q)} = O[\delta_n \sqrt{\zeta_n(4)}]$. So (59) follows from the identity

(61) $$N_n(\mathbf{g}) - N_n(0) = \sum_{I \subseteq \{1,\ldots,p\}} \int_0^{\mathbf{g}_I} \frac{\partial^{|I|} N_n(\mathbf{u}_I)}{\partial \mathbf{u}_I} d\mathbf{u}_I,$$

where the summation is over all the $2^p - 1$ nonempty subsets of $\{1, \ldots, p\}$. □

PROOF OF COROLLARY 1. The sequence $(\delta_n)$ clearly satisfies (16). Theorem 1 implies that $|\hat{\theta}_n| = O_\mathbb{P}(1) = o_\mathbb{P}(\delta_n)$. Note that $K_n(0) = \sum_{i=1}^n \psi(e_i)\mathbf{z}_{i,n}$ and $K_n(\hat{\theta}_n) = -\sum_{i=1}^n \varphi(-\mathbf{z}'_{i,n}\hat{\theta}_n)\mathbf{z}_{i,n} + O_\mathbb{P}(r_n)$. By Theorem 2, (19) follows from

$$\sum_{i=1}^n [\varphi(-\mathbf{z}'_{i,n}\hat{\theta}_n) + \mathbf{z}'_{i,n}\hat{\theta}_n\varphi'(0)]\mathbf{z}_{i,n} = \sum_{i=1}^n O[|\mathbf{z}'_{i,n}\hat{\theta}_n|^2]|\mathbf{z}_{i,n}| = O_\mathbb{P}[\zeta_n(3)]$$

in view of $\sum_{i=1}^n \mathbf{z}_{i,n}\mathbf{z}'_{i,n} = \mathrm{Id}_p$, $\zeta_n(2) = p$, $\zeta_n(3) \leq r_n\zeta_n(2) = O(r_n)$ and $\zeta_n(4) = O(r_n^2)$. For (20), it suffices to show that the left-hand side has the bound $O_\mathbb{P}[c_n\sqrt{\zeta_n(2+2\lambda)}\log n + c_nr_n]$ for any sequence $c_n \to \infty$. The latter easily follows from (19) and $m(t) = O(|t|^\lambda)$. □

4.1. *Proof of Theorem* 3.

LEMMA 6. *Under the assumptions of Theorem* 3, *we have* (i)

$$\sup_{|\beta| \leq b_n} |\tilde{K}_n(\beta) - \tilde{K}_n(0)| = O_{\text{a.s.}}(L_{\bar{n}} + B_{\bar{n}}) \tag{62}$$

*and* (ii) *for any* $\iota > 0$, $\tilde{K}_n(0) = o_{\text{a.s.}}(h_n)$, *where* $h_n = n^{1/2}(\log n)^{3/2} \times (\log \log n)^{1/2+\iota/4}$.

PROOF. As in the proof of Theorem 2, write $\tilde{K}_n = \tilde{M}_n + \tilde{N}_n$, where

$$\tilde{M}_n(\beta) = \sum_{i=1}^{n} \{\psi(e_i - \mathbf{x}_i'\beta) - \mathbb{E}[\psi(e_i - \mathbf{x}_i'\beta)|\mathcal{F}_{i-1}]\}\mathbf{x}_i,$$

$$\tilde{N}_n(\beta) = \sum_{i=1}^{n} \{\psi_1(-\mathbf{x}_i'\beta; \mathcal{F}_{i-1}) - \varphi(-\mathbf{x}_i'\beta)\}\mathbf{x}_i.$$

Since $n^{-5/2} = o(B_{\bar{n}})$, by Lemmas 7 and 8, (i) holds. For (ii), as with the argument in (42), we have $\|\tilde{K}_n(0)\| = O(\xi_n^{1/2}) = O(\sqrt{n})$. So the stated almost sure bound follows from the Borel–Cantelli lemma and (68) in view of the argument in (69). □

LEMMA 7. *Let* $(\pi_i)_{i \geq 1}$ *be a sequence of bounded positive numbers for which there exists a constant* $c_0 \geq 1$ *such that*

$$\max_{n \leq i \leq 2n} \pi_i \leq c_0 \min_{n \leq i \leq 2n} \pi_i \quad \text{holds for all large } n. \tag{63}$$

*Assume* (A5) *and* $\tilde{r}_n = O(\sqrt{n})$. *Let* $\varpi_d = 2c_0 \pi_{2^d}$ *and* $q > 3/2$. *Then as* $d \to \infty$,

$$\sup_{|\beta| \leq \varpi_d} \max_{n \leq 2^d} |\tilde{M}_n(\beta) - \tilde{M}_n(0)| = O_{\text{a.s.}}[\sqrt{\tilde{\tau}_{2^d}(\varpi_d)}d^q + 2^{-5d/2}]. \tag{64}$$

Lemma 7 is proved in [64]. The argument is similar to the one used for Lemma 4.

LEMMA 8. *Let* $(\pi_i)_{i \geq 1}$ *be a positive sequence satisfying* (63) *and* $\pi_n = o[n^{-1/2}(\log n)^2]$; *let* $\varpi_d = 2c_0 \pi_{2^d}$. *Assume* (21) *and assume* (A6) *holds with* $l = 0, \ldots, p$. *Then we have* (i)

$$\left\| \sup_{|\mathbf{g}| \leq \pi_n} |\tilde{N}_n(\mathbf{g}) - \tilde{N}_n(0)| \right\| = O[\sqrt{\xi_n(4)}\pi_n] \tag{65}$$

*and* (ii) *as* $d \to \infty$, *for any* $\iota > 0$,

$$\max_{n \leq 2^d} \sup_{|\mathbf{g}| \leq \varpi_d} |\tilde{N}_n(\mathbf{g}) - \tilde{N}_n(0)|^2 = o_{\text{a.s.}}[\xi_{2^d}(4)\varpi_d^2 d^3(\log d)^{1+\iota}]. \tag{66}$$

PROOF. Let $Q_{n,j}(\beta) = \sum_{i=1}^n \psi_1(-\mathbf{x}_i'\beta; \mathcal{F}_{i-1})x_{ij}$, $1 \leq j \leq p$, and $S_n(\beta) = Q_{n,j}(\beta) - Q_{n,j}(0)$. Since $\pi_n = o[n^{-1/2}(\log n)^2]$, by (21), $\pi_n \tilde{r}_{\bar{n}} \to 0$. It is easily seen that the argument in the proof of Lemma 5 implies that there exists a positive constant $C < \infty$ such that

$$\mathbb{E}\left\{\sup_{|\beta| \leq \varpi_d} |S_n(\beta) - S_{n'}(\beta)|^2\right\} \leq C \sum_{q=1}^p \varpi_d^{2q} \sum_{i=n'+1}^n |\mathbf{x}_i|^{2+2q} \quad (67)$$

holds uniformly over $1 \leq n' < n \leq 2^d$. So (65) holds. We now show (66). Let $\sup_\beta$ denote $\sup_{|\beta| \leq \varpi_d}$. By the maximal inequality (see [64])

$$\left\|\max_{n \leq 2^d} \sup_{|\beta| \leq \varpi} |S_n(\beta)|\right\| \leq \sum_{i=0}^d \left[\sum_{m=1}^{2^{d-i}} \mathbb{E}\left\{\sup_{|\beta| \leq \varpi} |S_{2^i m}(\beta) - S_{2^i(m-1)}(\beta)|^2\right\}\right]^{1/2}, \quad (68)$$

since $\iota > 0$ and $\varpi_d^{2q} \xi_{2^d}(2+2q) = O(\varpi_d^2 \xi_{2^d}(4))$, (67) implies that

$$\sum_{d=2}^\infty \frac{\|\max_{n \leq 2^d} \sup_{|\beta| \leq \varpi} |S_n(\beta)|\|^2}{\xi_{2^d}(4)\varpi_d^2 d^3 (\log d)^{1+\iota}} = \sum_{d=2}^\infty \frac{O[(d+1)^2]}{d^3 (\log d)^{1+\iota}} < \infty. \quad (69)$$

By the Borel–Cantelli lemma, (66) holds in view of (63). □

PROOF OF THEOREM 3. By Lemma 6, we have (i) and

$$\sup_{|\beta| \leq b_n} |\tilde{K}_n(\beta)| = O_{\text{a.s.}}(w_n), \quad \text{where } w_n = L_{\bar{n}} + h_n + B_{\bar{n}} \quad (70)$$

and $b_n = n^{-1/2}(\log n)^{3/2}(\log\log n)^{1/2+\iota}$. Let $\Theta_n(\beta) = \sum_{i=1}^n [\rho(e_i - \mathbf{x}_i'\beta) - \rho(e_i)]$ and

$$A_n(\beta) = -\sum_{i=1}^n \int_0^1 \varphi(-\mathbf{x}_i'\beta t)\mathbf{x}_i'\beta\, dt.$$

Using $\rho(e_i) - \rho(e_i - \mathbf{x}_i'\beta) = \int_0^1 \psi(e_i - \mathbf{x}_i'\beta t)\mathbf{x}_i'\beta\, dt$, we have by (70) that

$$\sup_{|\beta| \leq b_n} |\Theta_n(\beta) - A_n(\beta)| = \sup_{|\beta| \leq b_n} \left|\int_0^1 \tilde{K}_n(\beta t)\beta\, dt\right| = O_{\text{a.s.}}(w_n b_n).$$

Let $\lambda_* = 2^{-1} \liminf_{n \to \infty} \lambda_n/n$. By (21), $\tilde{r}_n b_n = o(1)$. Then $\xi_n(3)b_n^3 = O(n\tilde{r}_n)b_n^3 = o(nb_n^2)$. Since $\varphi(\delta) = \varphi(0) + \varphi'(0)\delta + O(\delta^2)$, we have for all large $n$ that

$$\inf_{|\beta|=b_n} A_n(\beta) \geq \tfrac{1}{3}\varphi'(0)n\lambda_* b_n^2.$$

Since $m(t) = O(\sqrt{t})$ as $t \to 0$, $L_n = O[\xi_n^{1/2}(3)b_n^{1/2}(\log n)^q]$. Let $3/2 < q < 7/4$. Under the condition $\xi_n(2) = O(n)$, we have $\xi_n(2+\kappa) \leq \tilde{r}_n^\kappa \xi_n(2) = O(n\tilde{r}_n^\kappa)$, $\kappa > 0$. Elementary calculations show that, with (21), $h_n = o(nb_n)$, $L_{\bar n} = o(nb_n)$ and $B_{\bar n} = o(nb_n)$. Therefore, $w_n = o(nb_n)$ and consequently we have

$$\inf_{|\beta|=b_n} \Theta_n(\beta) \geq \inf_{|\beta|=b_n} A_n(\beta) - \sup_{|\beta|\leq b_n} |\Theta_n(\beta) - A_n(\beta)|$$

$$\geq \tfrac{1}{3}\varphi'(0)n\lambda_* b_n^2 + O_{\text{a.s.}}(w_n b_n) \geq \tfrac{1}{4}\varphi'(0)n\lambda_* b_n^2$$

almost surely. By the convexity of the function $\Theta_n(\cdot)$,

$$\left\{\inf_{|\beta|\geq b_n} \Theta_n(\beta) \geq \tfrac{1}{4}\varphi'(0)n\lambda_* b_n^2\right\} = \left\{\inf_{|\beta|=b_n} \Theta_n(\beta) \geq \tfrac{1}{4}\varphi'(0)n\lambda_* b_n^2\right\}.$$

Therefore the minimizer $\hat\beta_n$ satisfies $|\hat\beta_n| \leq b_n$ almost surely.

(iii) Let $|\beta| \leq b_n$. Since $b_n \tilde r_n \to 0$, by Taylor's expansion,

$$-\sum_{i=1}^n \varphi(-\mathbf{x}_i'\beta)\mathbf{x}_i = \sum_{i=1}^n [\varphi'(0)\mathbf{x}_i'\beta + O(|\mathbf{x}_i'\beta|^2)]\mathbf{x}_i = \varphi'(0)\Sigma_n\beta + O[\xi_n(3)b_n^2].$$

So (24) follows from (i) and (ii) in view of $\tilde\Omega_n(\hat\beta_n) = O_{\text{a.s.}}(\tilde r_n)$. □

PROOF OF COROLLARY 2. Clearly the condition $\sup_u f_1(u|\mathcal{F}_0) \leq C_0 < \infty$ implies (A5). The other two conditions $\varphi(t) = t\varphi'(0) + O(t^2)$ and $m(t) = O(\sqrt{t})$ as $t \to 0$ easily follow. Then we have (A1)–(A3) and (A6). By Theorem 3, it remains to show that $\tilde\Omega_n(\hat\beta_n) = O_{\text{a.s.}}(\tilde r_n)$. Observe that $\psi(x) = \alpha - \mathbf{1}_{x\leq 0}$ and $\rho_\alpha(x+\delta) - \rho_\alpha(x) = \delta\psi(x) + \delta^+ \mathbf{1}_{x=0}$. Let $\Theta_n(\beta) = \sum_{i=1}^n [\rho_\alpha(e_i - \mathbf{x}_i'\beta) - \rho_\alpha(e_i)]$ and $\mathbf{v} = \tilde\Omega_n(\hat\beta_n)$. Then

$$\lim_{\epsilon\downarrow 0}\frac{\Theta_n(\hat\beta_n + \epsilon\mathbf{v}) - \Theta_n(\hat\beta_n)}{\epsilon} = \sum_{i=1}^n (-\mathbf{x}_i'\mathbf{v})\psi(e_i - \mathbf{x}_i'\hat\beta_n) + \sum_{i=1}^n (-\mathbf{x}_i'\mathbf{v})^+ \mathbf{1}_{e_i=\mathbf{x}_i'\hat\beta_n}$$

$$= -(\tilde\Omega_n(\hat\beta_n))'\mathbf{v} + \sum_{i=1}^n (-\mathbf{x}_i'\mathbf{v})^+ \mathbf{1}_{e_i=\mathbf{x}_i'\hat\beta_n}.$$

Since $\Theta_n(\hat\beta_n + \epsilon\mathbf{v}) \geq \Theta_n(\hat\beta_n)$ and $(-\mathbf{x}_i'\mathbf{v})^+ \leq |\mathbf{x}_i'||\mathbf{v}|$, we have $|\Omega_n(\hat\beta_n)| \leq \sum_{i=1}^n |\mathbf{x}_i'| \mathbf{1}_{e_i=\mathbf{x}_i'\hat\beta_n}$. By Lemma 9 and Schwarz's inequality, $|\tilde\Omega_n(\hat\beta_n)| \leq (p+1)\tilde r_n$ almost surely. □

LEMMA 9. *Assume that $\sup_u |f_1(u|\mathcal{F}_0)| \leq C_0$ almost surely for some constant $C_0$. Then*

$$(71) \qquad \sup_\beta \sum_{i=1}^n \mathbf{1}_{e_i=\mathbf{x}_i'\beta}|\mathbf{x}_i| \leq (p+1)\tilde r_n \qquad \text{almost surely.}$$

PROOF. It suffices to show $\mathbb{P}(e_{i_1} = \mathbf{x}'_{i_1}\beta, \ldots, e_{i_{p+1}} = \mathbf{x}'_{i_{p+1}}\beta \text{ for some } \beta) = 0$ for all $1 \le i_1 < \cdots < i_{p+1}$. To this end, the argument in [4] is useful. Clearly we can find $u_1, \ldots, u_{p+1}$ with $u_1^2 + \cdots + u_{p+1}^2 \ne 0$ such that $u_1\mathbf{x}_{i_1} + \cdots + u_{p+1}\mathbf{x}_{i_{p+1}} = 0$. Without loss of generality let $u_{p+1} = 1$ and write $\eta = \sum_{j=1}^{p} u_j e_{i_j}$. Then $\mathbb{P}(e_{i_{p+1}} = -\eta | \mathcal{F}_{i_{p+1}-1}) = 0$ almost surely. So $\mathbb{P}(u_1 e_{i_1} + \eta = 0) = 0$, which completes the proof. $\square$

**Acknowledgments.** I am grateful to the referees and an Associate Editor for their many helpful comments. I would also like to thank Jan Mielniczuk for valuable suggestions.

DEPARTMENT OF STATISTICS
UNIVERSITY OF CHICAGO
5734 S. UNIVERSITY AVENUE
CHICAGO, ILLINOIS 60637
USA
E-MAIL: wbwu@galton.uchicago.edu